\documentclass[10pt]{article}
\usepackage[english]{babel}										
\usepackage[utf8]{inputenc}										
\usepackage[T1]{fontenc}										
\usepackage{amsmath,amsfonts,amssymb,amsthm,cancel,siunitx,
calculator,calc,mathtools,empheq,latexsym, kotex}
\usepackage{ulem}
\usepackage{multirow}
\usepackage{float}
\usepackage{graphicx} 
\graphicspath{ {./images/} }
\usepackage{subfig,epsfig,tikz,float}		            
\usepackage{booktabs,multicol,multirow,tabularx,array}          
\newtheorem{problem}{Problem}
\newtheorem{hyp}{Hypothesis}
\newtheorem{remark}{Remark}
\title{\normalsize\bf%
\uppercase{Self-supervised learning for a nonlinear inverse problem with a forward operator involving an unknown function arising in photoacoustic tomography}
}
\author{%
Gyeongha Hwang$^{1}$, Gihyeon Jeon$^{2*}$, Sunghwan Moon$^{3}$
}

\begin{document}

\date{}

\maketitle

\vspace{-0.5cm}

\begin{center}
{\footnotesize 
$^1$ Department of Mathematics, Yeungnam University, Gyeongsan 38541, Republic of Korea\\
$^{2}$ School of Mathematics, Kyungpook National University, Daegu 41566, Republic of Korea\\ 
$^3$ Department of Mathematics, Kyungpook National University, Daegu 41566, Republic of Korea\\
*Corresponding author: rydbr6709@knu.ac.kr}
\end{center}

\bigskip
\noindent
{\small{\bf ABSTRACT.}
In this article, we are concerned with a nonlinear inverse problem with a forward operator involving an unknown function. The problem arises in diverse applications and is challenging in the presence of an unknown function, which makes it ill-posed. Additionally, the nonlinear nature of the problem makes it difficult to use traditional methods, and thus, the study addresses a simplified version of the problem by either linearizing it or assuming knowledge of the unknown function. Here, we propose self-supervised learning to directly tackle a nonlinear inverse problem involving an unknown function. In particular, we focus on an inverse problem derived in photoacoustic tomograpy (PAT), which is a hybrid medical imaging with high resolution and contrast. PAT can be modeled based on the wave equation. The measured data provide the solution to an equation restricted to surface and initial pressure of an equation that contains biological information on the object of interest. The speed of a sound wave in the equation is unknown. Our goal is to determine the initial pressure and the speed of the sound wave simultaneously. Under a simple assumption that sound speed is a function of the initial pressure, the problem becomes a nonlinear inverse problem involving an unknown function. The experimental results demonstrate that the proposed framework performs successfully.

}

\section{Introduction}\label{sec:1}
The inverse problem finds the cause factor from observed data, which has applications in fields such as optics, radar, acoustics, communication theory, signal processing, medical imaging, computer vision, geophysics, oceanography, and astronomy because it tells us about what we cannot directly observe. The forward operator (the inverse of the inverse problem) can be modeled as a (non)linear system and often involves an unknown function. Due to the nature of the inverse problem, it is usually very hard to know the cause factor. For example, in medical imaging the cause factor is the human body section, and in seismology, we never know the structure of the earth’s interior.

In this article, we are concerned with a nonlinear inverse problem with a forward operator involving an unknown function. Our goal is to simultaneously find, from the measurements, the unknown function and the inverse operator. The problem is generally ill-posed because of the unknown function. Additionally, nonlinearity in the problem makes conventional methods difficult to use. To handle the problem, one may simplify it linearly or assume knowledge about the unknown function. Here, we propose a self-supervised framework to directly tackle a nonlinear inverse problem involving an unknown function. In particular, we address an inverse problem derived in photoacoustic tomography (PAT). Although our framework is proposed to solve a problem arising in PAT, it is generic and can be extended to handle any nonlinear inverse problem involving an unknown function.

The rest of this section presents an introduction to PAT. In Section \ref{sec:problem_formulation}, we formulate the inverse problem arising in PAT, which is nonlinear and also involves an unknown function. The structure and learning method of the proposed framework for the problem are described in Section \ref{sec:network_design}. Numerical simulation results in Section \ref{sec:numerical_simulations} demonstrate that the proposed framework performs successfully. 

\subsection{Photoacoustic tomography}\label{PAT}
PAT is hybrid medical imaging that combines the high contrast of optical imaging with the high spatial resolution of ultrasound images~\cite{Jiang2018photoacoustic, Xia2014photoacoustic, kuchment2013radon}. The physical basis of PAT is the photoacoustic effect discovered by Bell in 1881~\cite{Bell1881production}. In PAT, when a non-destructive testing target object absorbs a non-ionizing laser pulse, it thermally expands and emits acoustic waves. The emitted ultrasound contains biological information on the target object, and is measured by a detector placed around it. The internal image of the target object is reconstructed from the measured data. 
The advantage of PAT is that it is economical and less harmful because of non-ionizing radiation use~\cite{Steinberg2019photoacoustic}. 

The propagation of the emitted ultrasound $p(\mathbf x,t)$ can be described by the wave equation:
\begin{equation}\label{eq:wave equation}
\partial_{t}^{2} p(\mathbf{x},t) = c(\mathbf{x})^{2} \Delta_{\mathbf{x}} p(\mathbf{x},t) \text{ on } \mathbb{R}^{2} \times [0, \infty),
\end{equation}
with initial conditions
\begin{equation}\label{eq:initial condition}
p(\mathbf{x},0) = f(\mathbf{x}) \qquad \partial_{t}p(\mathbf{x},0) = 0 \text{ on } \mathbb{R}^{2}.
\end{equation}
Here, $c$ is the speed of the waves, and $f$ is the initial pressure, which contains biological information such as the location of cancer cells in a physically small amount of tissue.
It is a natural assumption that $f$ has compact support in the bounded domain, $\Omega$, and the detectors are located on the boundary of the domain, $\partial \Omega$.
Regarding the measurement procedure, the point-shaped detector measures the average pressure above $\partial \Omega$ where the detectors are located, and this average pressure is the value of a pressure wave, $p(\mathbf{x},t)$.
Therefore, one of the mathematical problems in PAT is reconstructing $f$ from the measured data, $p|_{\partial \Omega \times [0, \infty)}$, which implies obtaining an internal image of the target object. 

It is well-known that given initial pressure $f$ and speed $c$, the solution, $p$, is determined uniquely. We define the wave's forward operator, $\mathcal{W}$, as follows:
$$
\mathcal{W} : (f, c) \mapsto p|_{\partial \Omega \times [0, \infty)},\quad \text{i.e.,}\quad  \mathcal{W}(f, c)=p|_{\partial \Omega \times [0, \infty)}.
$$
The reconstruction problem for $f$ from $\mathcal{W}(f, c)$ is studied when speed $c$ is constant~\cite{Zangerl2019photoacoustic,xu2005universal}. Oksanen and Uhlmann~\cite{oksanen2014photoacoustic} and Stefanov and Uhlmann~\cite{stefanov2009photoacoustic} studied explicit reconstruction when the sound speed is known. If $c$ depends on space variable $\mathbf{x}$, the problem become much more difficult. A few researchers have studied the problem with a given variable sound speed~\cite{qian2011efficient, belhachmi2016direct, hristova2008reconstruction, moon2022}. Liu and Uhlmann figured out the sufficient conditions for recovering $f$ and $c$~\cite{liu2015determining}. 

Recently, the application of deep learning in medical imaging, including PAT, has been investigated extensively. The roles of deep learning in tomography include forward and inverse operator approximation, image reconstruction from sparse data, and artifact/noise removal from reconstructed images~\cite{antholzer2018photoacoustic, Ongie2020deep, Wang2020deep, Grohl2021deep, Yang2021review, Antholzer2019deep, Zhou2021photoacoustic}. 
There are also studies on limited-view data~\cite{Guan2020limited, Zhang2020new}. Shan et al. proposed an iterative optimization algorithm that reconstructs $f$ and $c$ simultaneously via supervised learning~\cite{Shan2019simultaneous}.
However, most work deals with linear inverse problems or inverse problems without involving an unknown function~\cite{hoop2021}.

Many studies on PAT with deep learning are based on supervised learning. Supervised learning exploits a collection of data that pairs boundary data and initial pressure.  In practical applications, it is difficult to obtain the initial pressure, because initial pressure represents the internal human body. Therefore, it is necessary to study a learning method exploiting boundary data only. One such method is self-supervised learning that exploits supervised signals generated from input data by leveraging their structure~\cite{shurrab2022self, jing2020self}.

\section{Problem Formulation}\label{sec:problem_formulation}
In this section, we formulate the problem precisely. For this, we make several assumptions. 
First, we assume $f$ has compact support, since the target object is finite. 
Secondly, $c$ is assumed to be a function of $f$, namely $c(\mathbf{x})^{2} = \Gamma(f(\mathbf{x}))$ for some function $\Gamma:[0,1]\to [0,\infty)$, because wave speed $c$ depends on the medium. Lastly, we assume that $\Gamma(0)$ and $\Gamma(1)$ are known: $\Gamma(0)=c_0$ and $\Gamma(1)=c_1$. The last assumption is reasonable, because $\Gamma(0)$ and $\Gamma(1)$ represent wave speeds in the air and the highest thermal expansion coefficient, respectively.
Then, Equation \eqref{eq:wave equation} is rewritten as
\begin{equation}\label{eq:main_eqn}
\partial_{t}^{2} p(\mathbf{x},t) = \Gamma(f(\mathbf{x})) \Delta_{\mathbf{x}} p(\mathbf{x},t) \text{ on } \mathbb{R}^{2} \times [0, \infty).
\end{equation}
We define $\mathcal{W}_{\Gamma}$ as $\mathcal{W}_{\Gamma}(f) = p|_{\partial \Omega \times [0, \infty)}$ where $p$ is the solution of \eqref{eq:main_eqn} with initial conditions \eqref{eq:initial condition}. Then, the inverse problem can be formulated by determining unknown $\Gamma$ and $f$ from the given $\mathcal{W}_{\Gamma}(f)$. 
However, this problem is ill-posed; 
for any $\Gamma^{\prime}$ satisfying 
$$
\left\{
\begin{array}{ll}
\Gamma = \Gamma^{\prime} \text{ on } Im(f)\\
\Gamma \neq \Gamma^{\prime} \text{ on } Dom(\Gamma) \setminus Im(f)
\end{array}
\right.,
$$
we have $\mathcal{W}_{\Gamma}(f) = \mathcal{W}_{\Gamma^{\prime}}(f)$. 
Hence, $\Gamma$ cannot be uniquely determined from $\mathcal{W}_{\Gamma}(f)$. There is also a possibility that $\Gamma_{1}$, $\Gamma_{2}$, $f_{1}$, and $f_{2}$ exist such that $\Gamma_{1} \neq \Gamma_{2}$, $f_{1} \neq f_{2}$, and $\mathcal{W}_{\Gamma_{1}}(f_{1}) = \mathcal{W}_{\Gamma_{2}}(f_{2})$. Instead, we consider the following inverse problem.

\begin{problem}\label{prob:IP1}
Given that the collection of boundary data $\mathcal{B}_\Gamma := \{\mathcal{W}_{\Gamma}(f) \, |\,\Gamma : [0,1] \rightarrow [0,\infty), \Gamma(0)=c_0, \Gamma(1)=c_1 \mbox{ and } f \in L^{2}(\mathbb{R}^{2})\text{ has compact support}\}$,
\begin{enumerate}
    \item determine unknown $\Gamma$ from $\mathcal{B}_\Gamma$, and
    \item for all $\mathcal{W}_{\Gamma}(f) \in \mathcal{B}_\Gamma$, determine $f$.
\end{enumerate}
\end{problem}

Then, the uniqueness statements for Problem \ref{prob:IP1} are as follows:
\begin{hyp}\label{hyp:hyp1}
If $\Gamma_1 \neq \Gamma_2$, then $\mathcal{B}_{\Gamma_1}\neq \mathcal{B}_{\Gamma_2}.$
\end{hyp}

\begin{hyp}\label{hyp:hyp2}
For a fixed $\Gamma$, if $f_{1} \neq f_{2}$, then
$
\mathcal{W}_{\Gamma}(f_{1}) \neq \mathcal{W}_{\Gamma}(f_{2}).
$
\end{hyp}
In this article, we solve Problem \ref{prob:IP1} under hypothesis \ref{hyp:hyp1} and \ref{hyp:hyp2}. The problem is difficult to solve for two reasons:

\begin{enumerate}
    \item Expression \eqref{eq:main_eqn} involves an unknown $\Gamma$.
    \item Expression \eqref{eq:main_eqn} is not linear.
\end{enumerate}

We are going to solve Problem \ref{prob:IP1} by exploiting a deep neural network (DNN). Since DNN can only handle finite data, we address the following inverse problem.

\begin{problem}\label{prob:IP2}
For the given $\{\mathcal{W}_{\Gamma}(f_{i})\, |\, \Gamma : [0,1] \rightarrow [0,\infty), \Gamma(0)=c_0, \Gamma(1)=c_1 \mbox{ and } f_{i} \in L^{2}(\mathbb{R}^{2})\text{ has compact support}, i = 1,\cdots,N\}$, determine $\Gamma$ and $\{f_{i} | i = 1,\cdots,N\}$.
\end{problem}


\section{Network Design}\label{sec:network_design}
\begin{figure}[H]
    \centering
    \includegraphics[width=\textwidth]{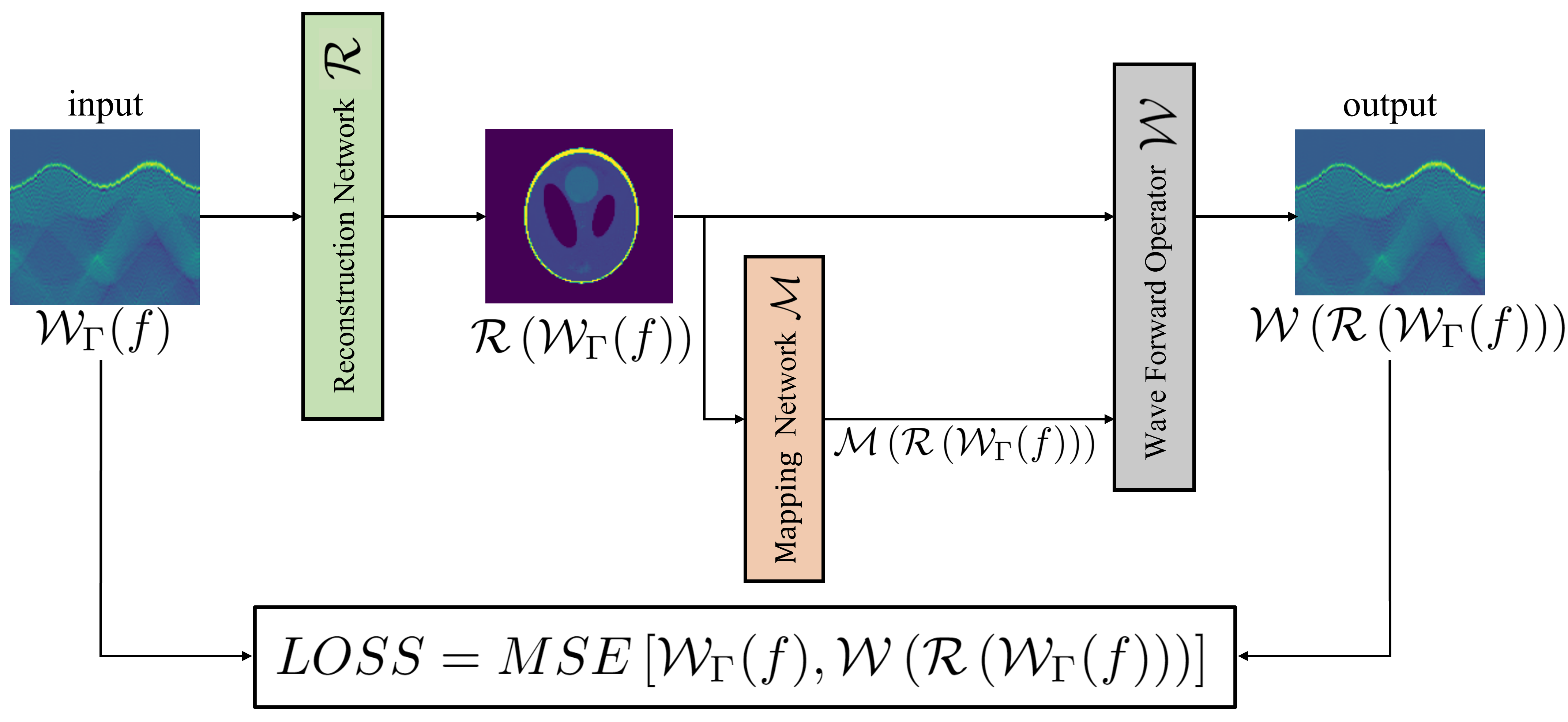}
    \caption{The proposed framework}
    \label{fig:network_design}
\end{figure}
We propose self-supervised learning for the problem formulated in Section \ref{sec:problem_formulation}. Our goal is simultaneously reconstructing $\{f_{i}\}_{i=1}^N$ and $\Gamma$ from given collection $\left\{\mathcal{W}_{\Gamma}(f_{i})\right\}_{i=1}^N$. The proposed framework is depicted in Figure \ref{fig:network_design}.
It consists of three components: 

\begin{enumerate}
    \item Reconstruction network $\mathcal R$
    \item Mapping network $\mathcal M$
    \item Wave forward operator $\mathcal W$.
\end{enumerate}

Reconstruction network $\mathcal R$ learns to reconstruct the initial data from the measured data. Mapping network $\mathcal M$ approximates the function $\Gamma : [0, 1] \to [0,\infty)$ satisfying $c(\mathbf{x})^{2} = \Gamma(f(\mathbf{x}))$. Forward operator $\mathcal W$ assigns the measured data to the initial data and the wave speed. Here, we adopt the $k$-space method. If every component in the framework functions properly, the output should be the same as the input. Thus, we define the loss function as the difference between input and output:
$$
\mathcal{L} = \dfrac{1}{N}\sum_{i=1}^N \dfrac{\left\| \mathcal{W}_{\Gamma}(f_i) - \mathcal{W}_{\mathcal M}(\mathcal R(\mathcal W_\Gamma (f_i)) \right\|_{2}}{\left\| \mathcal{W}_{\Gamma}(f_i) \right\|_{2}}.
$$
\begin{remark}
Our method estimates $\Gamma$ and the inverse operator, $\mathcal{W}_{\Gamma}^{-1}$. The estimated inverse operator can be used for fast inference of the initial pressure from the boundary measurement.
\end{remark}

\begin{remark}
The proposed framework is generic and can be extended to handle a nonlinear inverse problem involving an unknown function.
\end{remark}

Detailed structures of each component in the framework are described below.

\subsection{Reconstruction network $\mathcal R$}
\begin{figure}[H]
\includegraphics[width=\textwidth]{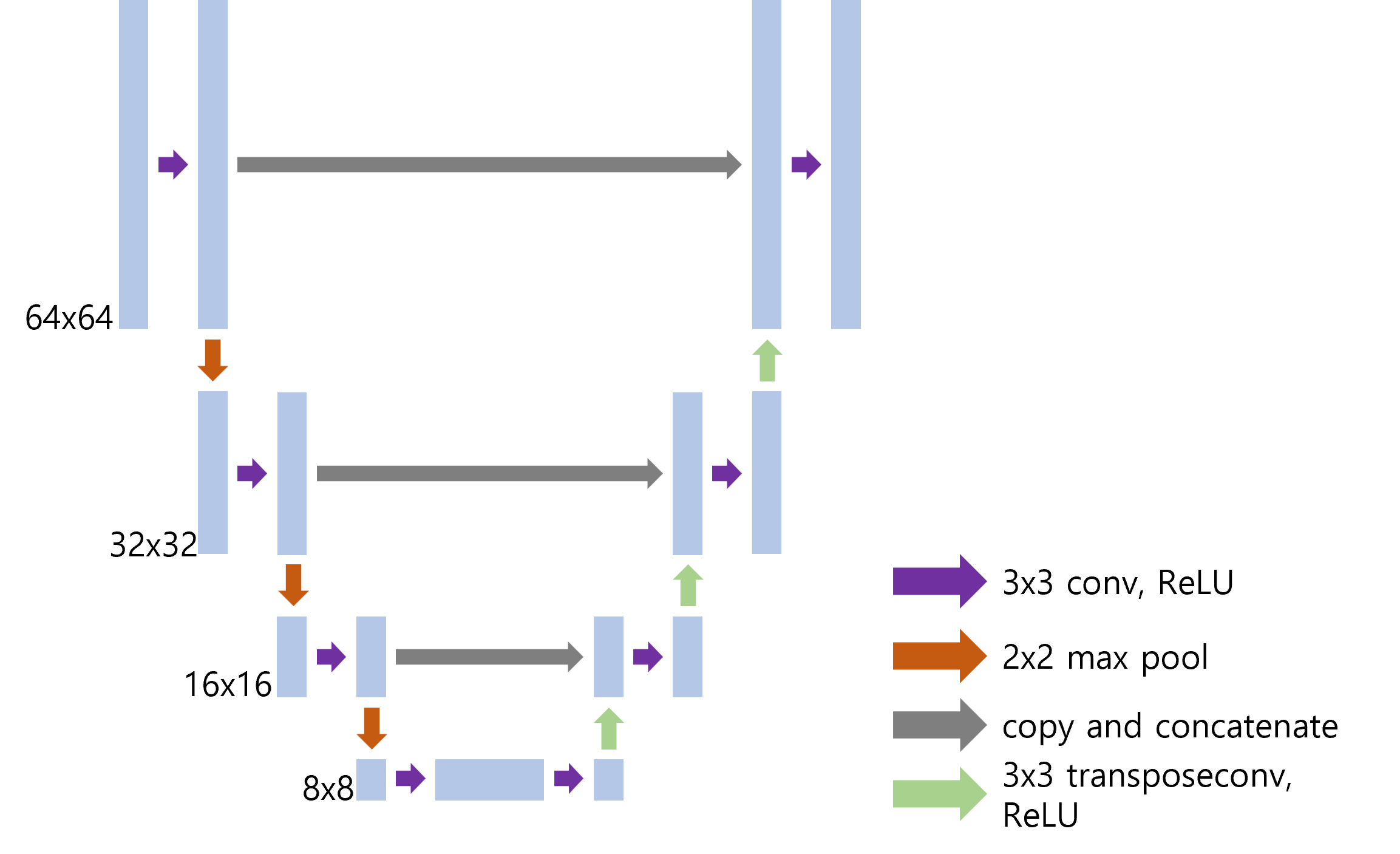}
\caption{The U-net architecture for data sized $64 \times 64$}
\label{fig:unet_design}
\end{figure}

Reconstruction network $\mathcal R$ reconstructs $f$ from input data $\mathcal{W}_{\Gamma}(f)$. Indeed, it approximates inverse map $\mathcal{W}_{\Gamma}^{-1} : \mathcal{W}_{\Gamma}(f) \mapsto f$. If speed $\Gamma$ of the wave is constant, it is well-known that the inverse map of Expression \eqref{eq:main_eqn} is linear~\cite{xu2005universal, moon2018inversion, gouia2021numerical}. Inspired by this fact, 
we propose the reconstruction network as a perturbation of a linear map:
\begin{equation}\label{eq:reconstruction_network}
\mathcal{R} := T_{1} + U \circ T_{2},    
\end{equation}
where  $T_{1}$, $T_{2} : \mathbb{R}^{m \times m} \to \mathbb{R}^{m \times m}$ are linear, and $U : \mathbb{R}^{m \times m} \to \mathbb{R}^{m \times m}$ is the U-net described in Figure \ref{fig:unet_design}. U-net is a type of convolutional neural network (CNN) introduced in~\cite{ronneberger2015u} and used widely in medical imaging. U-net consists of a contracting path and an expansive path. The contracting path has a typical CNN structure where the input data are extracted into a feature map with a small size and a large channel. In the expansive path, the size of the feature map increases again, and the number of channels decreases. At the end of $\mathcal{R}$, since the range of $f$ is $[0, 1]$, we use the clamp function, which rounds up values smaller than the minimum, and rounds down values larger than the maximum.

\begin{figure}[H]
\centering
\includegraphics[width=\textwidth]{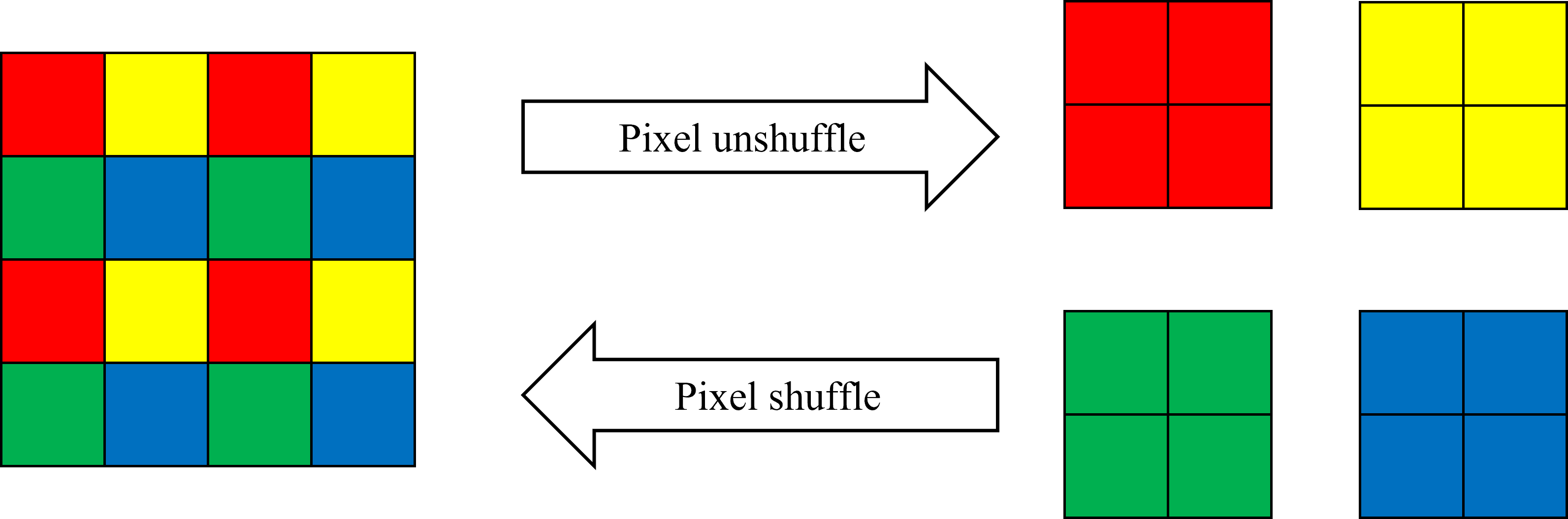}
\caption{Illustration of applying Pixel Shuffle and Pixel Unshuffle}
\label{fig:pixel}
\end{figure}

\begin{figure}[H]
\centering
\includegraphics[width=\textwidth]{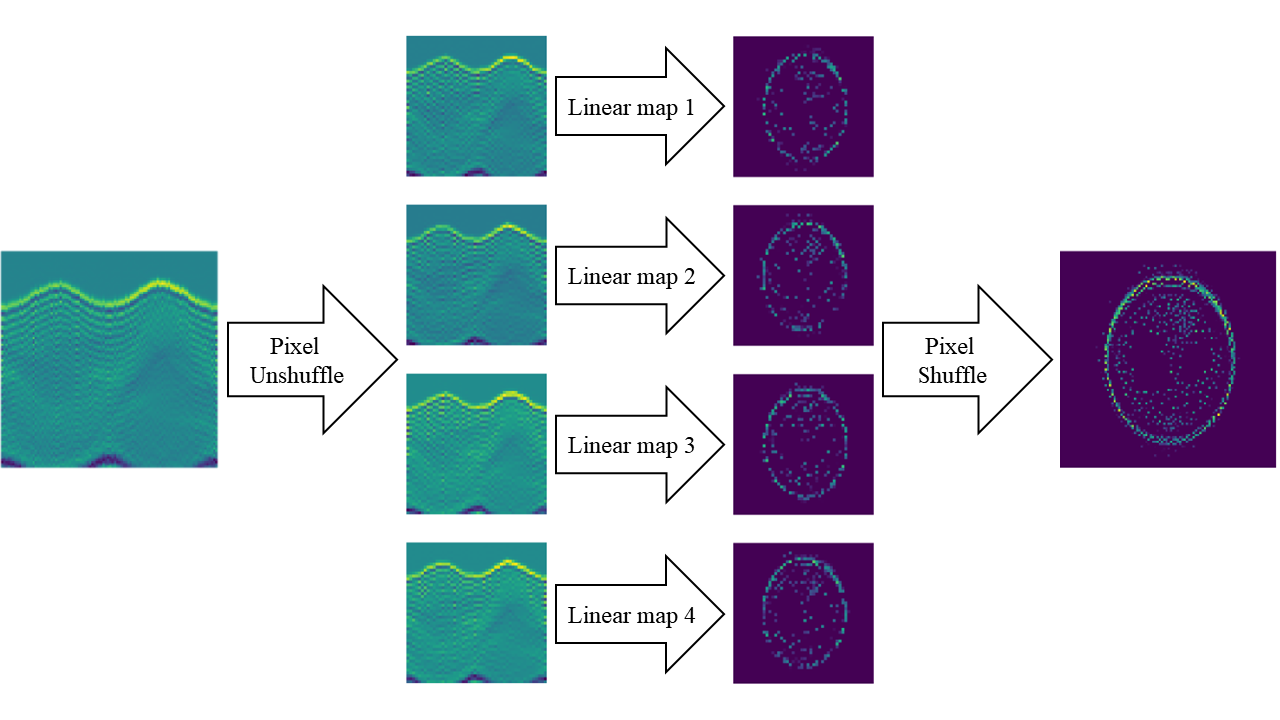}
\caption{Architecture of an alternative map of linear operators for high-resolution data}
\label{fig:alter}
\end{figure}

The proposed reconstruction network showed high performance with low-resolution data at $64 \times 64$ (see Section \ref{lowRes} below). With high-resolution data, however, the linear operators $T_{1}$ and $T_{2}$ in reconstruction network $\mathcal{R}$ create some problems because they contain too many parameters, causing a lot of critical points that impede convergence to the global minimum. They also create a hardware issue, and thus, for high-resolution data, we employ Pixel Shuffle and Pixel Unshuffle, which reduce the number of parameters contained in linear operators~\cite{shi2016real}. Pixel Unshuffle splits one image into several images, and Pixel Shuffle merges several images into one image, as illustrated in Figure \ref{fig:pixel}. Instead of applying the linear operators ($T_1$ and $T_2$) directly to high-resolution data, we process the data as follows (Figure \ref{fig:alter}) : 
\begin{enumerate}
\item Split high-resolution data ($m \times m$) into four sets of low-resolution data ($\frac{m}{2} \times \frac{m}{2}$) by exploiting Pixel Unshuffle.
\item Apply four different linear operators to the low-resolution data.
\item Merge the output of the linear operators by using Pixel Shuffle.
\end{enumerate}

\subsection{Mapping network $\mathcal{M}$}
We use multilayer perceptron (MLP) to approximate unknown $\Gamma$, because MLP can approximate any continuous function (for the universal approximation theorem, see \cite{cybenko1989approximation} and \cite{moon2022approximation}). 
The proposed network is a simple structure containing only three hidden layers of 10 nodes. To satisfy the assumption that $\Gamma(0)=c_{0}$ and $\Gamma(1)=c_{1}$, the output of MLP is slightly manipulated as follows:
$$
\mathcal{M}(f) = MLP(f) - MLP(0)*(1 - f) - MLP(1)*f + ((c_{1} - c_{0})f + c_{0}),
$$
so that
\begin{equation}\label{eq:constraints}
\mathcal{M}(0) = c_0 \text{ and }\mathcal{M}(1) = c_1.  
\end{equation}

\subsection{Forward problem}\label{sec:forward}
A solution to initial value problem $\eqref{eq:main_eqn}$ can be computed by the $k$-space method \cite{mast2001k,  cox2007k}. The $k$-space method is a numerical method for computing solutions to acoustic wave propagation, and it uses information in the frequency space to obtain a solution for the next time step. For calculating propagation of $p(\mathbf{x}, t)$, let $w(\mathbf{x}, t) = \dfrac{1}{\Gamma(f(\mathbf{x}))} p(\mathbf{x}, t)$ be an auxiliary field. Then, we have
$$
\partial_{t}^{2}w(\mathbf{x}, t) = \Delta_{\mathbf{x}} \left[ \Gamma(f(\mathbf{x})) w(\mathbf{x}, t) \right].
$$
Taking Fourier transform $\mathcal{F}_{\mathbf{x}}$ for $w$ with respect to $\mathbf{x}$ yields
\begin{equation}
\partial_{t}^{2} \mathcal{F}_{\mathbf{x}}w (\mathbf{k}, t) = -|\mathbf{k}|^{2}\mathcal{F}_{\mathbf{x}} \big[ \Gamma(f(\cdot)) w(\cdot, t) \big](\mathbf{k}).
\label{eq:fourier transform of W}
\end{equation}
Meanwhile, the numerical approximation of the second derivative of $\mathcal{F}_{\mathbf{x}}w$ is
\begin{equation}
\partial_{t}^{2}\mathcal{F}_{\mathbf{x}}w(\mathbf{k}, t) \approx \dfrac{\mathcal{F}_{\mathbf{x}}w(\mathbf{k}, t + \Delta t) - 2\mathcal{F}_{\mathbf{x}}w(\mathbf{k}, t) + \mathcal{F}_{\mathbf{x}}w(\mathbf{k}, t - \Delta t)}{(\Delta t)^{2}},
\label{eq:numerical approximation}
\end{equation}
where $\Delta t$ is the time step.
Then, by combining \eqref{eq:fourier transform of W} and \eqref{eq:numerical approximation}, we have
$$
\mathcal{F}_{\mathbf{x}}w(\mathbf{k}, t + \Delta t) = 2\mathcal{F}_{\mathbf{x}}w(\mathbf{k}, t) - \mathcal{F}_{\mathbf{x}}w(\mathbf{k}, t - \Delta t) - (\Delta t)^{2}
|\mathbf{k}|^{2}\mathcal{F}_{\mathbf{x}} \big[ \Gamma(f(\cdot)) w(\cdot, t) \big](\mathbf{k}).
$$
By taking the inverse Fourier transform, $\mathcal{F}^{-1}_{\mathbf{k}}$, we obtain
$$
w(\mathbf{x}, t + \Delta t) = 2w(\mathbf{x}, t) - w(\mathbf{x}, t - \Delta t) -  \mathcal{F}^{-1}_{\mathbf{k}} \bigg[ (\Delta t)^{2} |\mathbf{\cdot}|^{2}\mathcal{F}_{\mathbf{x}} \big[ \Gamma(f) w \big](\cdot, t) \bigg] (\mathbf{x}).
$$
Here, replacing $(\Delta t)^{2} |\mathbf{k}|^{2}$ in the third term with $4\sin^{2}\left(\frac{(\Delta t)|\mathbf{k}|}{2}\right)$ provides more accurate discretization \cite{mast2001k, cox2007k}. Finally, we have the wave propagation formula:
$$
w(\mathbf{x}, t + \Delta t) = 2w(\mathbf{x}, t) - w(\mathbf{x}, t - \Delta t) -  \mathcal{F}^{-1}_\mathbf{k} \bigg[ 4\sin^{2}\left(\frac{(\Delta t)|\cdot|}{2}\right)\mathcal{F}_\mathbf{x} \big[ \Gamma(f) w \big](\cdot, t) \bigg] (\mathbf{x}),
$$
or equivalently,
$$
p(\mathbf{x}, t + \Delta t)
= 2p(\mathbf{x}, t) - p(\mathbf{x}, t - \Delta t) -  \Gamma(f) \mathcal{F}^{-1}_\mathbf{k} \bigg[ 4\sin^{2}\left(\frac{(\Delta t)|\cdot|}{2}\right)\mathcal{F}_\mathbf{x} \big[ p \big](\cdot, t) \bigg] (\mathbf{x}).
$$

\section{Numerical Simulations}\label{sec:numerical_simulations}
In this section, we present the details from implementation of the proposed framework and from the experimental results when $\Omega$ is the unit ball.
\subsection{Datasets}
The Shepp-Logan phantom (an artificial image that describes a cross section of the brain) is commonly used for simulation in tomography and contains 10 ellipses \cite{shepp1974fourier}. Each ellipse is created with six parameters: the major axis, the minor axis, the $x$-coordinate and the $y$-coordinate of the center, the rotation angle, and the intensity value. 
The dataset of initial condition $f$ (defined on $[-1.0, 1.0]^{2} \subset \mathbb{R}^{2}$) is generated by slightly changing these six parameters with 
$$
\operatorname{supp}(f)\subset \left\{ (x_{}, y_{}) \in \mathbb{R}^{2} : \dfrac{{x_{}}^{2}}{0.69^{2}} + \dfrac{{y_{}}^{2}}{0.92^{2}} \le 1 \right\}.
$$ 
We created a set of 2,688 phantoms, $P = \{f_{i}\}_{i=1}^{2688}$. For $\Gamma$, we considered four cases: linear, square root, square, and constant:
\begin{enumerate}
    \item $\Gamma_{1}(f) = 0.3f + 0.7$
    \item $\Gamma_{2}(f) = 0.3\sqrt{f} + 0.7$
    \item $\Gamma_{3}(f) = 0.3f^{2} + 0.7$
    \item $\Gamma_{4}(f) = 0.7$
\end{enumerate} 
For $1 \le j \le 4$, we created a collection of data, $ \{\mathcal{W}_{\Gamma_{j}}f_{i}\}_{i=1}^{2688}$, by using the forward operator for $P$ and $\Gamma_{j}$. Of these data, we used 2,048 for training, 128 for validation, and 512 for testing.

\begin{figure}[H]
\includegraphics[width=\textwidth]{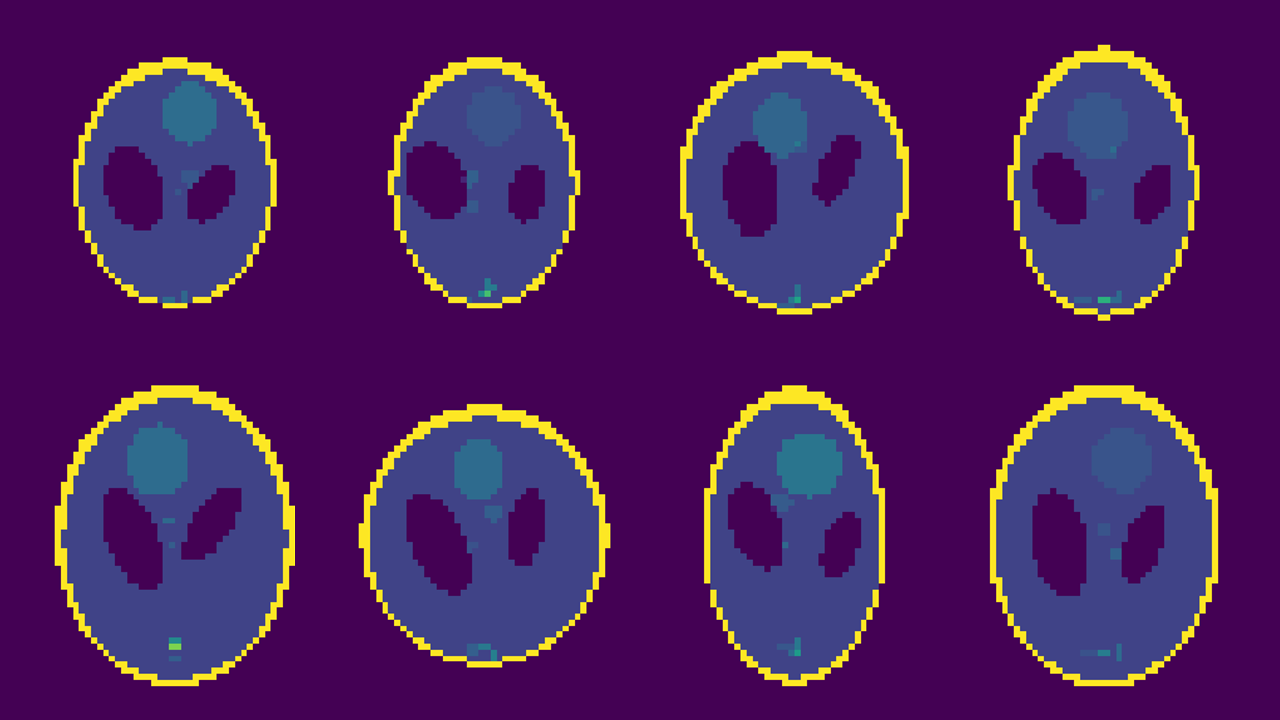}
\caption{Examples of phantoms}
\label{fig:datasets}
\end{figure}

\subsection{Training}

We used the Adam optimizer based on stochastic gradient descent and adaptive moment estimation to train the network \cite{kingma2014adam}. There are two neural networks in the proposed framework: reconstruction network $\mathcal{R}$ and mapping network $\mathcal{M}$. 
The learning rates for the linear term of $\mathcal{R}$, the perturbation term of $\mathcal{R}$, and for $\mathcal{M}$ were $10^{-4}$, $10^{-3}$, and $10^{-3}$, respectively. Momentum parameters of the Adam optimizer were set at $\beta_{1} = 0.9$ and $\beta_{2} = 0.999$.

We specifically set the batch size to 2. For general tasks, a moderately large batch size reduces the training time. However, in this problem, a small batch size is advantageous because our model must be able to reconstruct an exact image for the data, rather than an average result.

\subsection{Results}\label{sec:results}
In this section, we  illustrate the experimental results. The overall results are presented in Figure \ref{fig:mapping_plot_64}, Table \ref{table:error}, Figure \ref{fig:reconstruction_images_64}, Figure \ref{fig:mapping_plot_96}, Table \ref{table:error_for_96}, and Figure \ref{fig:reconstruction_images_96}. 
Here, the losses for $f$ and $\mathcal{W}_\Gamma(f)$ are respectively defined by
$$
\text{loss for }f =  \dfrac{1}{N}\sum_{i=1}^N \dfrac{\| f_i - \mathcal R(\mathcal W_\Gamma (f_i)) \|_{2}}{\| f_i\|_{2}},
$$
and 
$$
\text{loss for }\mathcal{W}_\Gamma{f} =  \dfrac{1}{N}\sum_{i=1}^N \dfrac{\left\| \mathcal{W}_{\Gamma}(f_i) - \mathcal{W}_{\mathcal M}(\mathcal R(\mathcal W_\Gamma (f_i)) \right\|_{2}}{\left\| \mathcal{W}_{\Gamma}(f_i) \right\|_{2}}.
$$

\subsubsection{Low-resolution data}\label{lowRes}
\begin{figure}[H]
\includegraphics[width=\textwidth]{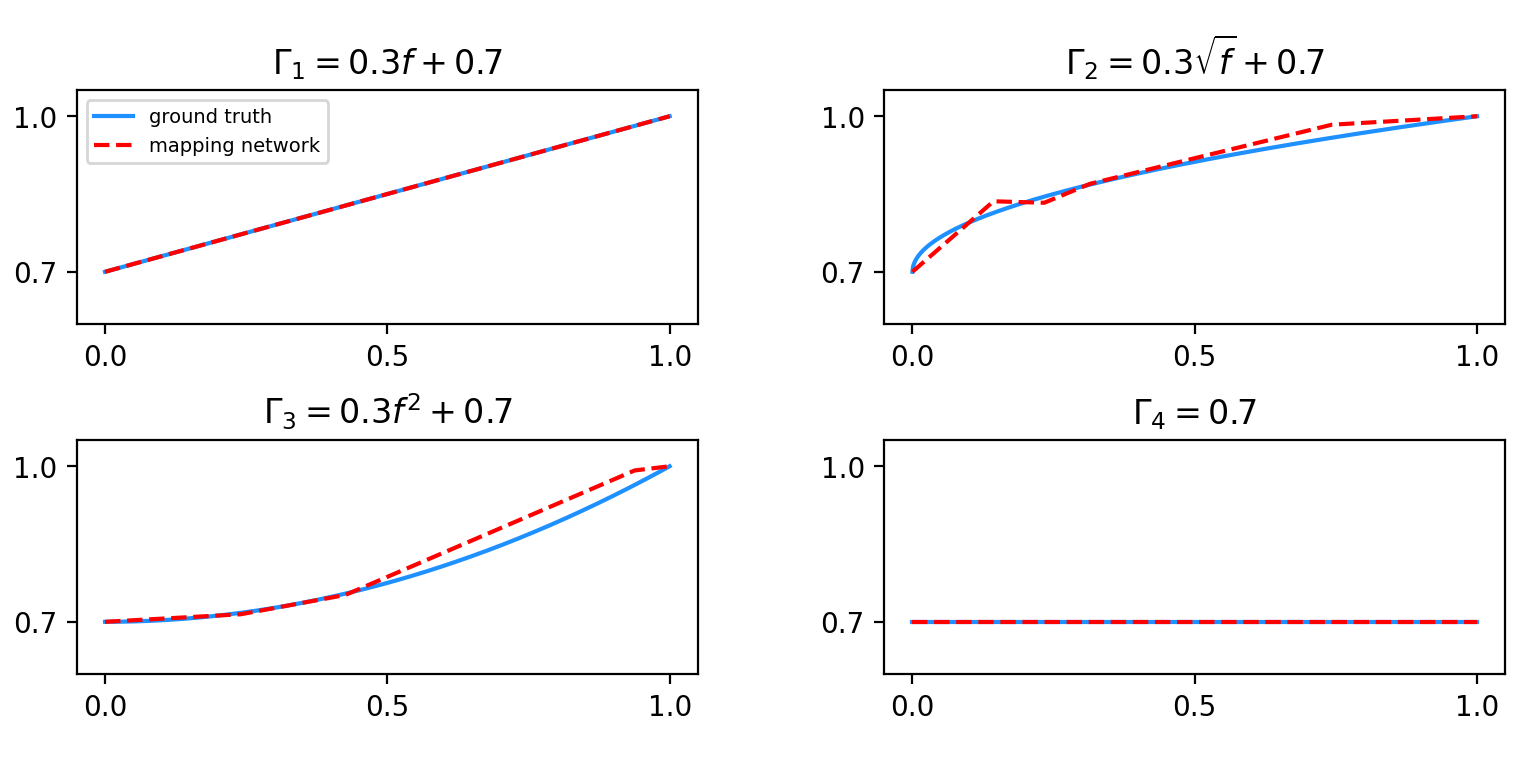}
\caption{Comparison of mapping network $\mathcal{M}$ and ground truth $\Gamma$ for data sized $64 \times 64$}
\label{fig:mapping_plot_64}
\end{figure}

\begin{figure}[H]
\includegraphics[width=\textwidth]{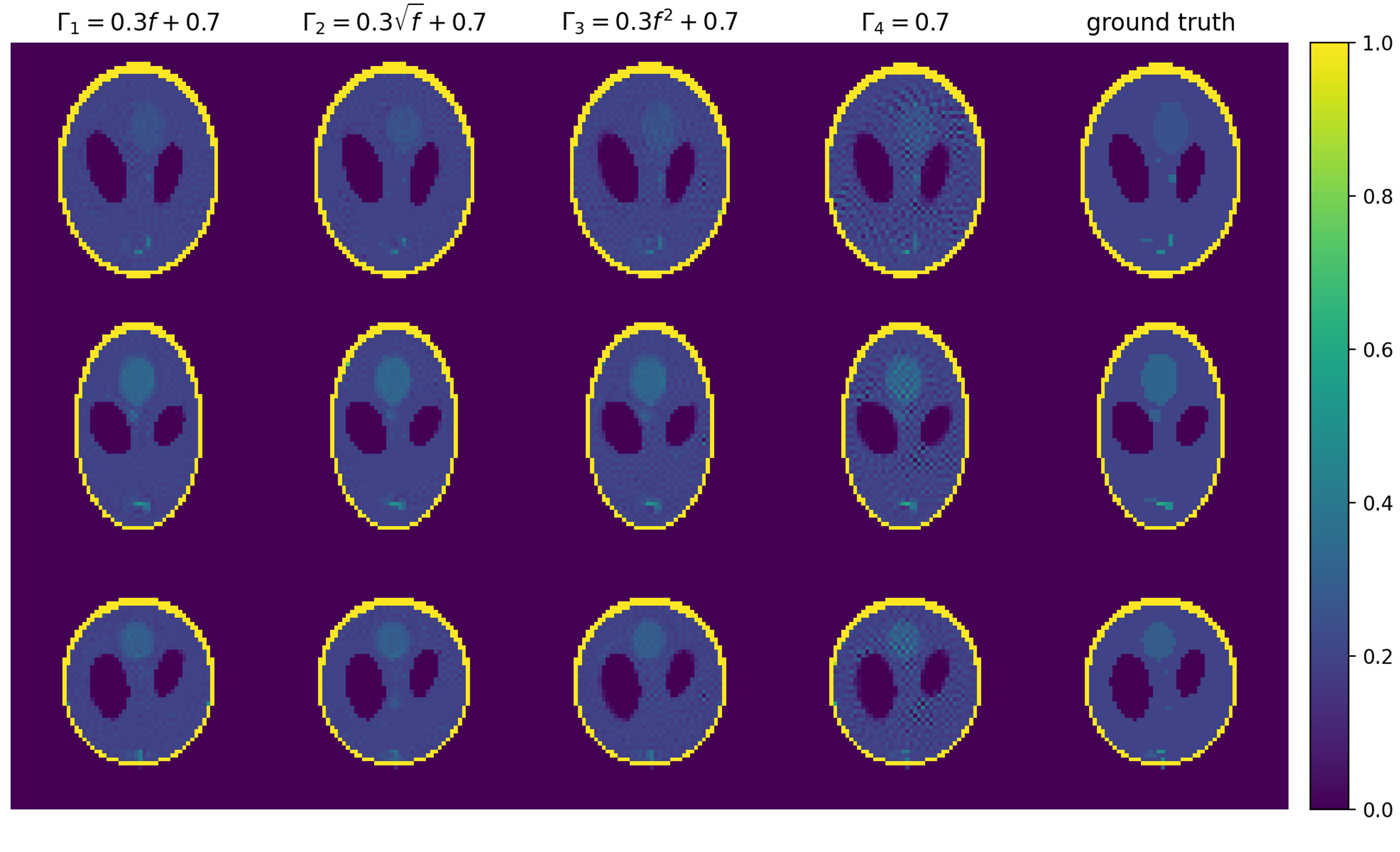}
\caption{Reconstruction results according to $\Gamma$ for data sized $64 \times 64$}
\label{fig:reconstruction_images_64}
\end{figure}

We conducted a simulation utilizing a dataset of images sized $64 \times 64$. Results from the mapping network are shown in Figure \ref{fig:mapping_plot_64}. We see that the mapping network accurately approximates $\Gamma$. When $\Gamma_{3} = 0.3f^{2} + 0.7$, there is a difference between the plot of mapping network $\mathcal{M}$ and the plot of $\Gamma$. This is because the values of $f \in P$ almost all belong to $[0, 0.3] \cup {1}$, so they have little effect on $\mathcal{W}_\Gamma(f)$. In all cases, the process of training the mapping network requires approximately $10^3$ iterations. Results from the reconstruction network are presented in Table \ref{table:error} and Figure \ref{fig:reconstruction_images_64}. Table \ref{table:error} shows the test errors. So we can conclude that the reconstruction network accurately approximates the inverse map in each case. Training the reconstruction network requires approximately $10^5$ iterations. 

\begin{remark}
The assumption for $\Gamma$ in \eqref{eq:constraints} is crucial. If constraint \eqref{eq:constraints} is not put on $\mathcal{M}$, it may take a long time to approximate $\Gamma$, or it may fail to find $\Gamma$. Under the constraint, $\mathcal{M}$ can quickly determine $\Gamma$. Early determination of $\Gamma$ helps to learn reconstruction network. 
\end{remark}

\renewcommand{\arraystretch}{1.5}
\begin{table}[ht]
\centering
\begin{tabular}[t]{lcc}
\hline
Assumption & loss for $f$ & loss for $\mathcal{W}_{\Gamma}(f) $\\
\hline
$\Gamma_{1} = 0.3f + 0.7$ & 0.00504 & 0.00702 \\
$\Gamma_{2} = 0.3\sqrt{f} + 0.7$ & 0.00537 & 0.00947 \\
$\Gamma_{3} = 0.3f^2 + 0.7$ & 0.00557 & 0.00634 \\
$\Gamma_{4} = 0.7$ & 0.01373 & 0.00456 \\
\hline
\end{tabular}
\caption{Test errors for $f$ and $\mathcal{W}_{\Gamma}(f)$ according to $\Gamma$ after 102,400 iterations for data sized $64 \times 64$}
\label{table:error}
\end{table}%

\subsubsection{High-resolution data}
\begin{figure}[H]
\includegraphics[width=\textwidth]{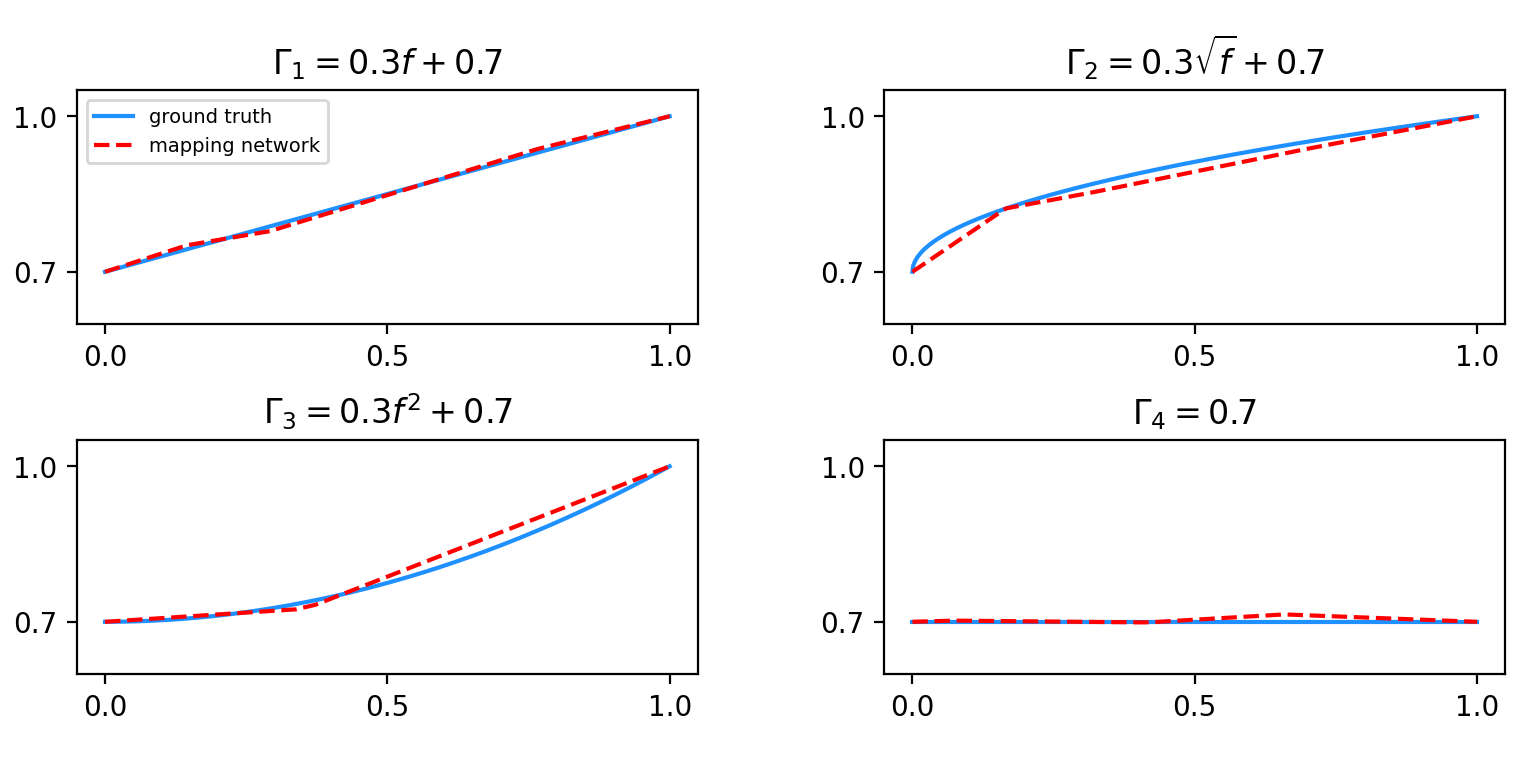}
\caption{Comparison of mapping network $\mathcal{M}$ and ground truth $\Gamma$ for data sized $96 \times 96$}
\label{fig:mapping_plot_96}
\end{figure}

\begin{figure}[H]
\includegraphics[width=\textwidth]{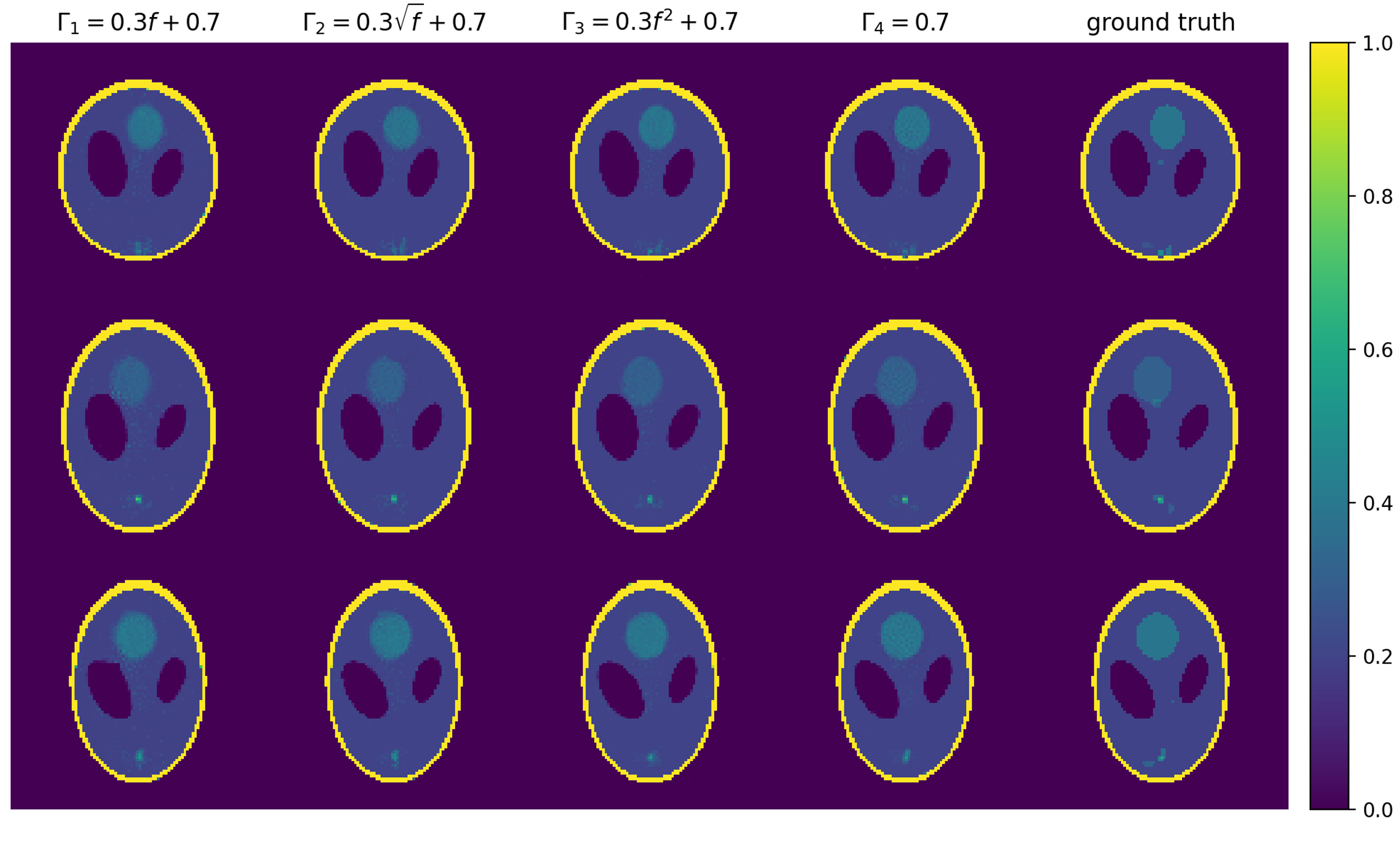}
\caption{Reconstruction results based on $\Gamma$ for data sized $96 \times 96$}
\label{fig:reconstruction_images_96}
\end{figure}

In the simulation for high-resolution data, two linear operators, $T_{1}$ and $T_{2}$ for the reconstruction network expressed in \eqref{eq:reconstruction_network}, are replaced by the alternative map described in Figure \ref{fig:alter}. The dataset was prepared with images sized at $96 \times 96$. Similar to the case with low-resolution data, mapping network $\mathcal{M}$ approximates $\Gamma$ accurately within $10^3$ iterations (Figure \ref{fig:mapping_plot_96}). On the other hand, for each $\Gamma$, the reconstruction network exhibits a slight decrease in performance that is acceptable (Table \ref{table:error_for_96} and Figure \ref{fig:reconstruction_images_96}). We surmise that the slight decrease in performance result from the reduction in parameters from applying Pixel Unshuffle and Pixel Shuffle.

\renewcommand{\arraystretch}{1.5}
\begin{table}[ht]
\centering
\begin{tabular}[t]{lcc}
\hline
Assumption & loss for $f$ & loss for $\mathcal{W}_{\Gamma}(f) $\\
\hline
$\Gamma_{1} = 0.3f + 0.7$ & 0.00860 & 0.01293 \\
$\Gamma_{2} = 0.3\sqrt{f} + 0.7$ & 0.01023 & 0.01679 \\
$\Gamma_{3} = 0.3f^2 + 0.7$ & 0.00710 & 0.01132 \\
$\Gamma_{4} = 0.7$ & 0.00689 & 0.69511 \\
\hline
\end{tabular}
\caption{Test errors for $f$ and $\mathcal{W}_{\Gamma}(f)$ according to $\Gamma$ after 102,400 iterations for data sized $96 \times 96$}
\label{table:error_for_96}
\end{table}%

\section{Conclusions}\label{sec:4}
We proposed self-supervised learning for a nonlinear inverse problem with a forward operator involving an unknown function. In medical imaging such as PAT, the initial pressure is mostly untrackable for the measured data. Moreover, it is difficult to know the wave speed. So, it is necessary to reconstruct initial pressure $f$ and the wave speed simultaneously. Under the simple assumption, the problem becomes a nonlinear inverse problem involving an unknown function. The experimental results demonstrated high performance from the proposed framework, which can be extended to a nonlinear inverse problem involving an unknown function and formulated under more complicated situations. This can be an interesting line of future research.

\section{Acknowledgement}
G. Hwang was supported by the 2019 Yeungnam University research grant. The work of G. Jeon and S. Moon was supported by the National Research Foundation of Korea (NRF-2022R1C1C1003464).

\bibliographystyle{unsrt}

\end{document}